\documentclass[letterpaper, 10 pt, journal, twoside]{IEEEtran}  

\IEEEoverridecommandlockouts                              

\usepackage{graphics} 
\usepackage{amsmath} 
\usepackage{amssymb}  
\usepackage{theorem}
\usepackage{color,xspace}
\usepackage{hyperref}
\usepackage{graphicx}
\usepackage{cite}
\usepackage{cleveref}
\usepackage{bm}
\usepackage{bbm}
\usepackage{tikz}

\newtheorem{theorem}{Theorem}
\newtheorem{lemma}{Lemma}	

\newtheorem{remark}{Remark}

\newtheorem{definition}{Definition}

\newcommand*\bmat[1]{\begin{bmatrix}#1\end{bmatrix}}

\def\dom{\mathrm{dom}}

\def\R{\mathbb{R}}
\newcommand{\prox}{\mathrm{prox}}
\title{A Control-Theoretic Approach to Analysis and Parameter Selection of Douglas-Rachford Splitting}

\author{Jacob H. Seidman$^1$, Mahyar Fazlyab$^2$, Victor M. Preciado$^2$, George J. Pappas$^2$
\thanks{1: Department of Applied Mathematics and Computational Science, University of Pennsylvania Email: seidj@sas.upenn.edu. 2: Department of Electrical and Systems Engineering, University of Pennsylvania. Email: \{mahyarfa, preciado, pappasg\}@seas.upenn.edu. Work was supported by ARL CRA DCIST W911NF-17-2- 0181 and the NSF under grants CAREER-ECCS-1651433.}} 

\begin{document}

\maketitle
\thispagestyle{empty}
\pagestyle{empty}

\begin{abstract}

Douglas-Rachford splitting and its equivalent dual formulation ADMM are widely used iterative methods in composite optimization problems arising in control and machine learning applications. The performance of these algorithms depends on the choice of step size parameters, for which the optimal values are known in some specific cases, and otherwise are set heuristically.  We provide a new unified method of convergence analysis and parameter selection by interpreting the algorithm as a linear dynamical system with nonlinear feedback.  This approach allows us to derive a dimensionally independent matrix inequality whose feasibility is sufficient for the algorithm to converge at a specified rate.  By analyzing this inequality, we are able to give performance guarantees and parameter settings of the algorithm under a variety of assumptions regarding the convexity and smoothness of the objective function. In particular, our framework enables us to obtain a new and simple proof of the $O(1/k)$ convergence rate of the algorithm when the objective function is not strongly convex.
\end{abstract}
\section{Introduction}
\label{Introduction}

\IEEEPARstart{I}{n this} paper, we consider problems of the form
\begin{align} \label{eq: main prob}
\mathrm{minimize}_{x \in \mathbb{R}^d} \ \{F(x)=f(x)+g(x)\},
\end{align}
where $f,g:\R^d \to \R \cup \{+\infty\}$ are convex, closed, and proper (c.c.p.).  Douglas-Rachford splitting (DRS) solves problem \eqref{eq: main prob} with the following iterations:
\begin{subequations} \label{eq: DR algorithm}
\begin{align}
y_k &= \prox_{\alpha f}(x_k),  \label{eq: proximal f update} \\
z_k &=\prox_{\alpha g}(2y_k-x_k),  \label{eq: proximal g update}  \\
x_{k+1} &= x_k + \lambda_k (z_k-y_k),  \label{eq: x update} 
\end{align}
\end{subequations}
where $\prox$ is the proximal operator (see Definition \ref{def: prox}) and $\alpha$ and $\lambda_k$ are known as the proximal step size and relaxation parameter, respectively.  For a proper selection of these parameters, the limiting values of both $y_k$ and $z_k$ will be a solution to \eqref{eq: main prob}.  The goal of this work is to provide convergence rates for DRS over various assumptions on $f$ and $g$, and optimize these rates with respect to the algorithm parameters $\alpha$ and $\lambda_k$ using semidefinite programs (SDPs).

The algorithm was first proposed in \cite{douglas1956numerical}, and has since found application in general separable optimization problems \cite{stathopoulos2016operator}.  Its dual formulation, ADMM, has been particularly useful in distributed optimization problems \cite{boyd2011distributed}.  Since the iterates of ADMM can be written as applying DRS to the dual problem \cite{gabay1983chapter, eckstein1989splitting}, convergence results for one algorithm are valid for the other as well when strong duality holds.

The convergence of DRS has previously been analyzed using monotone operator theory and variational inequalities, see \cite{lions1979splitting}, \cite{eckstein1992douglas}.  These techniques have led to proofs of a $O(1/k)$ convergence rate for the non-strongly convex case \cite{he20121}, and linear convergence when $f$ is smooth and strongly convex \cite{deng2016global}, \cite{giselsson2017linear}.  In the most general case, a condition for convergence is that $\lambda_k \in (0,2)$ with $\sum_{k=0}^\infty \lambda_k(2-\lambda_k) = \infty$ \cite{combettes2004solving}, though there exist cases where the algorithm converges with $\lambda_k > 2$ \cite{nishihara2015general}.

Recently, there has been interest in automating the analysis and design of optimization algorithms via SDPs, \cite{drori2014performance,lessard2016analysis,fazlyab2018analysis,hu2017dissipativity,van2018fastest,ryu2018operator,fazlyab2018design}. In particular, 
%
through the method of integral quadratic constraints proposed in \cite{lessard2016analysis}, the authors of \cite{nishihara2015general} derive an SDP for choosing the parameters of ADMM in the case of smooth and strongly convex $f$. Using a similar framework, the authors of \cite{francca2016explicit} provide evidence that as the relaxation parameter approaches $2$ from below, the linear convergence rate is close to being optimal and in \cite{francca2017tuning} are able to analytically solve the SDP to give a convergence rate.  The work in \cite{ghadimi2015optimal} gives an optimal choice for the relaxation parameter when $f$ is quadratic.  Furthermore, \cite{giselsson2014diagonal} gives a set of assumptions in which a bound on the linear convergence rate is minimized by setting $\lambda_k = 2$. 

\emph{Our Contribution:} By viewing DRS as a linear system with non-linear feedback, we derive a dimensionally independent matrix inequality which gives convergence guarantees via Lyapunov functions.  Whereas such an approach was previously applied in \cite{nishihara2015general} to the case of smooth and strongly convex $f$, our framework is novel in that it encompasses varying assumptions on the smoothness and convexity of $f$.  By changing a single term in the Lyapunov function for each scenario, we are able to relate the satisfaction of a matrix inequality to the convergence of the algorithm.  In particular, we give a new and simple proof of $O(1/k)$ convergence in the non-strongly convex case.  These symbolic results can then be used to select step sizes that optimize the derived rates.

In the strongly convex case, the corresponding matrix inequality is sufficient to guarantee a linear convergence rate.  We are able to modify the matrix inequality to linearize the dependence on $\lambda_k$, allowing us to numerically optimize its value for the convergence rate directly.  While previous work derived SDP's which can verify the performance of the algorithm for a given parameter setting, to the best of our knowledge this is the first time such a method immediately gives an optimal relaxation parameter when solved numerically, as opposed to having to search over a range of values for $\lambda_k$.

\section{Preliminaries}

We denote the set of real numbers by $\mathbb{R}$, the set of real $n$-dimensional vectors by $\mathbb{R}^n$, the set of real $m\times n$-dimensional matrices by $\mathbb{R}^{m\times n}$, and the $n$-dimensional identity matrix and zero matrix by $I_n$ and $0_n$, respectively.  For a function $f \colon \mathbb{R}^d \to \mathbb{R} \cup \{+\infty\}$, we denote by $\dom \, f=\{x \in \mathbb{R}^d \colon f(x)<\infty\}$ the effective domain of $f$. The subdifferential of a function $f$ at a point $x$ is $\partial f(x) := \{g\;|\;f(y) - f(x) \geq g^\top(y - x), \forall y \in \dom(f)\}$.  By abuse of notation we will also refer to a subgradient, that is an element of the subdifferential by $\partial f(x)$ as well.  The indicator function of a set $C$ is given by $\mathrm{1}_C(x) = 0$ if $x \in C$ and $\mathrm{1}_C(x) = \infty$ if $x \notin C$.  For two matrices $A \in \mathbb{R}^{m\times n}$ and $B\in \mathbb{R}^{p\times q}$ their Kronecker product is $A \otimes B$.  
 
We say a differentiable function $f \colon \mathbb{R}^d \to \mathbb{R}$ is $L_f$-smooth on $\mathcal{S} \subseteq \dom \, f$ if $\|\nabla f(x)-\nabla f(y)\|_2 \leq L_f \|x-y\|_2$ for some $L_f>0$ and all $x,y \in \mathcal{S}$. This also implies for all $x,y \in \mathcal{S}$, $f(y) \leq f(x) + \nabla f(x)^\top (y-x) +(L_f/2) \|y-x\|_2^2$.  A differentiable function $f \colon \mathbb{R}^d \to \mathbb{R}$ is $m_f$-strongly convex on $\mathcal{S} \subseteq \dom \, f$ if $m_f \|x-y\|_2^2 \leq (x-y)^\top (\nabla f(x)-\nabla f(y))$ for some $m_f>0$ and all $x,y \in \mathcal{S}$.  The class of functions which are $L_f$-smooth and $m_f$-strongly convex is denoted by $\mathcal{F}(m_f,L_f)$. 

\begin{definition}[Proximal Operator] \normalfont \label{def: prox}
	Given a c.c.p. function $f:\R^d \to \R \cup \{+\infty\}$ and $\alpha > 0$, the proximal operator $\prox_{\alpha f}:\R^d \to \R^d$ is defined as
	\begin{align}
	\prox_{\alpha f}(x) = \text{argmin}_y \left\{f(y) + \frac{1}{2\alpha}\|x - y\|_2^2\right\}.
	\end{align}
	The point $y = \prox_{\alpha f}(x)$ also is given by the implicit solution to the subgradient equation
	\begin{align} \label{eq: prox_implicit}
	y = x - \alpha \partial f(y).
	\end{align}
\end{definition}

We say that a nonlinear function $\phi:\R^d \to \R^d$ satisfies the incremental quadratic constraint \cite{accikmecse2011observers} (or point-wise integral quadratic constraint \cite{lessard2016analysis}) defined by $Q\in \R^{2d \times 2d}$ if for all $x,y$,
\begin{align} \label{eq: QC def}
	\begin{bmatrix}
	x\!-\!y \\ \phi(x)\!-\!\phi(y)
	\end{bmatrix}^\top \! Q \! \begin{bmatrix}
	x\!-\!y \\ \phi(x)\!-\!\phi(y)
	\end{bmatrix} \geq 0.
\end{align}
%
For $0\leq m, L < \infty$ define
	\begin{align} \label{eq: quadratic matrix}
	Q(m,L) = \begin{bmatrix}
	-\frac{m L}{m+L} & 1/2 \\ 1/2 & -\frac{1}{m+L}
	\end{bmatrix} \otimes I_d,
	\end{align}
and define $Q(m,\infty)$ as $\lim_{L \to \infty} Q(m,L)$.  It was noted in \cite{nesterov2013introductory,lessard2016analysis} that a differentiable function $f$ belongs to the class $\mathcal{F}(m_f,L_f)$ on $\mathcal{S}$ if and only if the gradient $\nabla f$ satisfies the incremental QC defined by $Q(m_f,L_f)$.  When $L_f = \infty$, the subgradient $\partial f$ satisfies the QC defined by $Q(m_f,\infty)$.  If we define
\begin{align}\label{eq: prox QC}
	Q_p(m,L,\alpha) = \begin{bmatrix}
	0  &  I_d \\ \alpha I_d & -I_d
	\end{bmatrix}Q(m,L) \begin{bmatrix}
	0  &  \alpha I_d \\ I_d & -I_d
	\end{bmatrix},
	\end{align}
then the proximal operator of a function $f \in \mathcal F(m,L)$, $\prox_{\alpha f}$, satisfies the incremental QC defined by $Q_p(m,L,\alpha)$ \cite{fazlyab2018analysis}.


\section{Analysis of Douglas-Rachford Splitting via Matrix Inequalities}

\subsection{Douglas-Rachford Splitting as a Dynamical System}
We can write the updates in \eqref{eq: DR algorithm} as a linear system with state $x_k$ and feedback nonlinearity $\phi(x_k)$,
\begin{align}
x_{k+1} = x_k + \lambda_k \phi(x_k),
\end{align}
where
\begin{align}
\phi(x_k) := \prox_{\alpha g}(2\prox_{\alpha f}(x_k)-x_k)-\prox_{\alpha g}(x_k).
\end{align}

Our main technique is to describe the nonlinearity $\phi$ with incremental QCs representing the $\prox$ operator.  This allows us to derive a matrix inequality as a sufficient condition for closed-loop stability of the system via a Lyapunov function argument.  We perform this derivation in the following three cases:
\begin{itemize}
\item \emph{Case 1:} $f \in \mathcal F(0,\infty)$ and $g \in \mathcal F(0,\infty)$,

\item \emph{Case 2:} $f \in \mathcal F(0,L_f)$ and $g \in \mathcal F(0,\infty)$, with $0 < L_f < \infty$,

\item \emph{Case 3:} $f \in \mathcal F(m_f,L_f)$ and $g \in \mathcal F(0,\infty)$, with $0 < m_f \leq L_f < \infty$.
\end{itemize}

We will see that for each case only one term in the Lyapunov function needs to be modified to obtain the convergence result.  We then use the matrix inequality condition for each case to obtain information about optimal choices of the algorithm parameters both symbolically and numerically.

\subsection{Characterization of Fixed Points}

From relation \eqref{eq: prox_implicit}, the iterates \eqref{eq: DR algorithm} can be rewritten as
\begin{subequations} \label{eq: implicit_updates}
\begin{align}
y_k &= x_k -\alpha \partial f(y_k),  \\
z_k &=2y_k-x_k - \alpha \partial g(z_k), \\
x_{k+1} &= x_k + \lambda_k (z_k-y_k) .
\end{align}
\end{subequations}
The fixed points of \eqref{eq: implicit_updates} satisfy
\begin{align} \label{eq: fixed points of DR}
x_{\star} = y_{\star} \!+\! \alpha \partial f(y_{\star}), \ y_{\star}=z_{\star}, \ \partial f(y_{\star})\! +\! \partial g(z_{\star})=0.
\end{align}
Since $y_\star = z_\star$, the rightmost equality is exactly the optimality condition for \eqref{eq: main prob}. 

We will also make use of the following relation, obtained from adding the equations in \eqref{eq: implicit_updates} and the definition of $\phi$,
\begin{align}\label{eq: zy partial relation}
\phi(x_k) = z_k - y_k = -\alpha(\partial f(y_k) + \partial g(z_k)).
\end{align}
From this, we can interpret the feedback nonlinearity $\phi$ as the optimality residual of problem \eqref{eq: main prob}, which is driven to zero by the linear system in the feedback interconnection.
%
\subsection{Convergence Certificates via Matrix Inequalities}

\subsubsection{Case 1: Non-strongly convex and non-smooth case}

We first assume that $f,g \in \mathcal{F}(0,\infty)$.  We propose the following family of Lyapunov functions parameterized by a sequence $\{\theta_i\}_{i=0}^\infty$ with $\theta_i > 0$,
%
\begin{align} \label{eq: Lyapunov}
V_k \!=\! \|x_k-x_{\star}\|_2^2+ \sum_{i=0}^{k-1} \theta_i \|\partial f(y_i) + \partial g(z_i)\|_2^2, 
\end{align}
for all $k > 0$ and $V_0 = \|x_0 - x_{\star}\|^2$.  For notational convenience define the partial sums $\Theta_k = \sum_{i=0}^{k-1} \theta_i$.  The presence of the running sum of subgradients is reminiscent of the Popov criterion \cite{haddad2011nonlinear}.  It can also be interpreted as the running weighted sum of fixed point residuals (see \eqref{eq: fixed points of DR}).  The next lemma shows how this Lyapunov function can ensure a convergence rate in terms of the growth of $\Theta_k$.
\begin{lemma} \label{lemma: weak lyapunov}
	Consider the algorithm in \eqref{eq: DR algorithm}. Suppose there exists a sequence $\{\theta_i\}_{i=0}^\infty$ with $\theta_i>0$ such that $V_{k+1} \leq V_k$ for all $k\geq 0$. Then 
	\begin{align} \label{eq: weak lyapunov rate}
\min_{i=0,\ldots,k-1} \|\partial f(y_i) + \partial g(z_i)\|_2^2 \leq \frac{1}{ \Theta_k}\|x_0-x_{\star}\|_2^2.
\end{align}
	
\end{lemma}
\begin{IEEEproof}
Since $V_{k+1} \leq V_k$ for all $k$, in particular we have that $V_k \leq V_0$, or
\begin{align}
\|x_k\!-\!x_{\star}\|_2^2\!+\!  \! \sum_{i=0}^{k-1} \theta_i\|\partial f(y_i) \!+\! \partial g(z_i)\|_2^2 \leq \|x_0-x_{\star}\|_2^2.
\end{align}
Removing the first term on the left, and dividing through by $\Theta_k$ gives
\begin{align} 
	 \sum_{i=0}^{k-1}\frac{\theta_i}{\Theta_k} \|  \partial f(y_i) + \partial g(z_i)\|_2^2 \leq \dfrac{\|x_0-x_ \star\|_2^2}{\Theta_k}.
	\end{align}
	The result follows from the fact that the left side is a weighted average, as $\theta_i > 0$ and $\sum_{i=0}^{k-1}\theta_i/\Theta_k = 1$.
\end{IEEEproof}

In the following theorem, we derive a matrix inequality in terms of $\alpha,\lambda$, and $\{\theta_i\}_{i=0}^\infty$ as a sufficient condition to guarantee $V_{k+1} \leq V_k$, which in turn implies \eqref{eq: weak lyapunov rate}.
%
%
\begin{theorem}\label{thm: weak thm}
	Let $m_f = 0$, $L_f = \infty$, and consider the following matrix inequality
	\begin{align}
	W_{k}^{(0)} + \sigma_k^{(1)} Q^{(1)} + \sigma_k^{(2)} Q^{(2)}  \preceq 0,\label{eq: weak LMI}
	\end{align}
	where 
\begin{subequations}\label{eq: Q1 and Q2}
\begin{align} 
W_k^{(0)} &= \begin{bmatrix}
0 & -\lambda_k & \lambda_k \\ -\lambda_k & \lambda_k^2 +\frac{\theta_k}{\alpha^2} & -\left(\lambda_k^2 + \frac{\theta_k}{\alpha^2}\right) \\ \lambda_k  & -\left(\lambda_k^2 +\frac{\theta_k}{\alpha^2}\right) & \lambda_k^2 +\frac{\theta_k}{\alpha^2}
\end{bmatrix} \!\otimes \!I_d \label{eq: W_0}, \\
Q^{(1)} &= \begin{bmatrix} 0 & I_d \\  \alpha I_d & - I_d \\ 0 & 0
\end{bmatrix}Q(m_f,L_f)  \begin{bmatrix} 0 & \alpha I_d & 0 \\  I_d & - I_d & 0
\end{bmatrix}, \\
Q^{(2)} &=\begin{bmatrix} 0 & -I_d \\ 0 & 2 I_d \\  \alpha I_d & - I_d
\end{bmatrix}\!Q(0,\infty) \!
\begin{bmatrix} 0 \!& 0\! &  \alpha I_d\! \\ - I_d \!& 2 I_d\! & - I_d \! \end{bmatrix}.
\end{align}
\end{subequations}
	If  $\alpha, \lambda_k,\theta_k > 0$ and $\sigma_k^{(1)}, \sigma_k^{(2)} \geq 0$ are chosen so that \eqref{eq: weak LMI} is satisfied for all $k \geq 0$, then for all $f \in \mathcal{F}(0,\infty)$ and $g \in \mathcal{F}(0,\infty)$ the iterates in \eqref{eq: DR algorithm} satisfy 
\begin{align}
\min_{i=0,\ldots,k-1} \|\partial f(y_i) + \partial g(z_i)\|_2^2 \leq \frac{1}{\Theta_k}\|x_0-x_{\star}\|_2^2.
\end{align}
\end{theorem}
\begin{IEEEproof}
We first see that $V_{k+1} - V_k$ can be written as a quadratic form.  Define the error signal
\begin{align} \label{eq: error signal}
e_k :=\begin{bmatrix}
(x_k-x_{\star})^\top & (y_k-y_{\star})^\top & (z_k-z_{\star})^\top \end{bmatrix}^\top.
\end{align}
Using the updates in \eqref{eq: implicit_updates}, the fact that $z_\star = y_\star$ (see \eqref{eq: fixed points of DR}), and the relation \eqref{eq: zy partial relation}, it can be verified that
\begin{align}
V_{k+1} - V_k = e_k^\top W_k^{(0)} e_k,
\end{align}
where $W_k^{(0)}$ is given by \eqref{eq: W_0}. Next, note that
\begin{align}
e_k^\top Q^{(1)} e_k = \bmat{x_k - x_\star \\ y_k - y_\star}^\top Q_p(m_f,L_f,\alpha) \bmat{x_k - x_\star \\ y_k - y_\star},
\end{align}  
where $Q_p(m_f,L_f,\alpha)$ is defined in \eqref{eq: prox QC}.  Since $y_k = \prox_{\alpha f}(x_k)$ and $y_\star =\prox_{\alpha f}(x_\star)$, this is exactly the incremental QC that the $\prox_{\alpha f}$ operator satisfies.  Thus, we have for all $k$, $e_k^\top Q^{(1)} e_k \geq 0$. We also note that
\begin{align}
\begin{bmatrix} 0\! &\! 0 &\!  \alpha I_d\! \\ - I_d\! & 2 I_d\! & - I_d  \!\end{bmatrix}e_k \! = \!\bmat{0 \!& \!\alpha \!I_d \!\\ I_d\!& -I_d\!}\!\bmat{(2y_k \!-\! x_k) \!-\! (2y_\star\! -\!x_\star) \!\\ z_k - z_\star }.\nonumber
\end{align}

As $z_k = \prox_{\alpha g}(2y_k - x_k)$ and $z_\star = \prox_{\alpha g}(2y_\star - x_\star)$, we similarly conclude that $e_k^\top Q^{(2)} e_k \geq 0$ is implied from the incremental QC that $\prox_{\alpha g}$ satisfies. Returning to \eqref{eq: weak LMI}, if we multiply from the left and right by $e_k^\top$ and $e_k$ respectively, we obtain
\begin{align}
e_k^\top W_k^{(0)} e_k + \sigma_k^{(1)} e_k^\top Q^{(1)} e_k + \sigma_k^{(2)} e_k^\top Q^{(2)} e_k \leq 0.
\end{align}
Since $\sigma_k^{(1)}, \sigma_k^{(2)} \geq 0$ and we have shown that $e_k^\top Q^{(1)} e_k \geq 0$ and $e_k^\top Q^{(2)} e_k \geq 0$, it must be that $e_k^\top W_k^{(0)} e_k \leq 0$.  Hence, $V_{k+1} \!-\! V_k \leq 0$, and the result now follows from Lemma \ref{lemma: weak lyapunov}.
\end{IEEEproof}

\subsubsection{Case 2: Non-strongly convex and smooth $f$}
If $f \in \mathcal F(0,L_f)$ with $0< L_f< \infty$, we may leverage the smoothness of $f$ to refine the result of the previous section.  In particular, we can use the inequality for $L_f$-smooth functions (see Preliminaries) to relate the behavior of the subgradients to the objective values, whereas in the previous section this inequality was not available.  For $k > 0$ let 
\begin{align} \label{eq: weak smooth lyapunov}
V_k = \|x_k-x_{\star}\|_2^2 +  \sum_{i=0}^{k-1} \theta_i[F(z_i)-F(z_{\star})],
\end{align}
with $V_0$ defined as in the previous case.

This Lyapunov function leads to the following Lemma, the proof of which is identical to that of Lemma \ref{lemma: weak lyapunov}.

\begin{lemma} \label{lemma: weak smooth lyapunov}
	Consider the algorithm in \eqref{eq: DR algorithm}. Suppose there exists a sequence $\{\theta_i\}_{i=0}^\infty$ with $\theta_i>0$ such that $V_{k+1} \leq V_k$ for all $k \geq0$. Then 
	\begin{align*}
	\min_{i=0,\ldots,k-1} [F(z_i) - F(z_\star)] \leq \dfrac{1}{\Theta_k} \|x_0 - x_\star\|_2^2.
	\end{align*}
\end{lemma}

This allows us to prove the following theorem for when $f$ is non-strongly convex and smooth.

\begin{theorem} \label{thm: weak smooth thm}
Let $0 = m_f < L_f < \infty$ and consider the following matrix inequality
	\begin{align}
	W_k^{(1)} + \sigma_k^{(1)} Q^{(1)} + \sigma_k^{(2)} Q^{(2)}  \preceq 0,\label{eq: weak smooth LMI}
	\end{align}
	where 
	\begin{align} \label {eq: W_1}
W_k^{(1)}\!\! =\! \!\begin{bmatrix} 0\!\! & -\lambda_k & \lambda_k \\
	-\lambda_k\!\! & \frac{\theta_k L_f}{2} + \lambda_k^2 & \frac{\theta_k}{2}\left(\frac{1}{\alpha} - L_f\right) \!-\!\lambda_k^2\\
	\lambda_k\!\! & \frac{\theta_k}{2}\!\left(\frac{1}{\alpha} \!- \!L_f\right)\! -\!\lambda_k^2 &  \theta_k\left(\frac{L_f}{2} - \frac{1}{\alpha}\right)\! +\!\lambda_k^2\end{bmatrix} \!\otimes\! I_d,
\end{align}
	and $Q^{(1)}$ and $Q^{(2)}$ are defined in \eqref{eq: Q1 and Q2}.  If $\alpha, \lambda_k,\theta_k > 0$ and $\sigma_k^{(1)}, \sigma_k^{(2)} \geq 0$ are chosen so that \eqref{eq: weak smooth LMI} is satisfied for all $k \geq 0$, then for all $f \in \mathcal{F}(0,L_f)$ with $0<L_f<\infty$ and $g \in \mathcal{F}(0,\infty)$ the iterates in \eqref{eq: DR algorithm} satisfy
\begin{align} \label{eq: weak smooth rate}
\min_{i=0,\ldots,k-1} [F(z_i) - F(z_\star)] \leq \frac{1}{\Theta_k} \|x_0-x_{\star}\|_2^2.
\end{align}
\end{theorem}

\begin{IEEEproof}
We begin by bounding the difference of the Lyapunov function defined in \eqref{eq: weak smooth lyapunov}, $V_{k+1} - V_k$, by a quadratic form in the error signal $e_k$ (see \eqref{eq: error signal}).   From the convexity and smoothness of $f$, we can write
\begin{align}
f(z_k) \!- \!f(y_k)\! &\leq \nabla f(y_k)^\top(z_k \!-\! y_k) +\frac{L_f}{2}\|z_k - y_k\|_2^2, \\
f(y_k) \!- \!f(z_\star)\! &\leq \nabla f(y_k)^\top(y_k \!-\!y_\star),
\end{align}
where we have used that $z_\star = y_\star$.  From the convexity of $g$
\begin{align}
g(z_k) -g(z_\star) \leq \partial g(z_k)^\top(z_k - z_\star).
\end{align}
Adding these three inequalities together and using the relation \eqref{eq: zy partial relation} allows us to conclude
\begin{align}\label{eq: weak smooth bound}
F(z_k) - F(z_\star) &\leq \frac{L_f}{2}\|y_k-y_\star\|_2^2 + \left(\frac{L_f}{2} - \frac{1}{\alpha}\right)\|z_k - z_\star\|_2^2\nonumber\\
&\quad+ \left(\frac{1}{\alpha} - L_f\right)(y_k-y_\star)^\top(z_k - z_\star).
\end{align}
Using the recursion for $x_{k+1}$, we then find that 
\begin{align} \label{eq: lyapunov W_1 rel}
V_{k+1} - V_k \leq e_k^\top W_k^{(1)} e_k.
\end{align}

The proof now proceeds identically as in the proof of Theorem \ref{thm: weak thm} up to the statement that \eqref{eq: weak smooth LMI} implies $e_k^\top W_k^{(1)} e_k \leq 0$.  Then by \eqref{eq: lyapunov W_1 rel}, we have that $V_{k+1} - V_k \leq 0$, and the result follows from Lemma \ref{lemma: weak smooth lyapunov}.
\end{IEEEproof}

\subsubsection{Case 3: Strongly convex and smooth $f$}
We now assume that $f \in \mathcal F(m_f,L_f)$ and $g \in \mathcal F(0, \infty)$, with $0 < m_f \leq L_f < \infty$.  For this scenario let
%
\begin{align} \label{eq: strong_lyapunov}
V_k = \|x_k - x_\ast\|_2^2.
\end{align}
The following lemma characterizes when we can extract a linear convergence rate from this Lyapunov function.
\begin{lemma} \label{lemma: lin lyapunov}
Consider the algorithm in \eqref{eq: DR algorithm}. Suppose there exists $ \rho \in (0,1)$ such the Lyapunov function $V_k$ defined by \eqref{eq: strong_lyapunov} satisfies $V_{k+1} \leq \rho^2 V_k$ for all $k \geq 0$.  Then 
\begin{align} \label{eq: linear_rate}
\|x_k-x_\ast \|_2^2\leq  \rho^{2k} \|x_0 - x_\ast\|_2^2.
\end{align}
\end{lemma}
\begin{IEEEproof}
The proof follows immediately $V_{k+1} \leq \rho^2 V_k$, the definition \eqref{eq: strong_lyapunov} of $V_k$, and induction.
\end{IEEEproof}

We again see that the difference $V_{k+1} - \rho^2 V_k$ can be written as a quadratic form acting on the error signal $e_k$ as defined in \eqref{eq: error signal}.  Using the definition for $x_{k+1}$ in terms of the previous iterates, we can write
\begin{align} \label{eq: strong_lyapunov_quadratic_relation}
V_{k+1} - \rho^2 V_k = e_k^\top Q_ke_k,
\end{align}
where $Q_k$ is given by
\begin{align}
Q_k = \begin{bmatrix}
1- \rho^2 & -\lambda_k & \lambda_k \\ -\lambda_k & \lambda_k^2 & -\lambda_k^2 \\ \lambda_k  & -\lambda_k^2& \lambda_k^2
\end{bmatrix} \otimes I_d. \label{eq: Q def}
\end{align}
From \eqref{eq: strong_lyapunov_quadratic_relation}, we can use the same reasoning developed in the non-strongly convex case to arrive at the following theorem.

\begin{theorem}\label{thm: strong thm}
	Let $0 < m_f \leq L_f < \infty$ and consider the following matrix inequality
	\begin{align} 
	Q_k + \sigma_k^{(1)} Q^{(1)} + \sigma_k^{(2)} Q^{(2)}  \preceq 0,\label{eq: strong LMI}
	\end{align}
	where $Q_k$ is given in \eqref{eq: Q def} and $Q^{(1)}$ and $Q^{(2)}$ are given in \eqref{eq: Q1 and Q2}.  If $\alpha, \lambda_k > 0$, $\sigma_k^{(1)}, \sigma_k^{(2)} \geq 0$, and $ \rho \in (0,1)$ are chosen so that \eqref{eq: strong LMI} is satisfied for all $k \geq 0$, then for all $f \in \mathcal{F}(m_f,L_f)$ and $g \in \mathcal{F}(0,\infty)$ with $0<m_f\leq L_f<\infty$, the iterates in \eqref{eq: DR algorithm} satisfy the following linear convergence rate 
\begin{align}
\|x_k - x_\ast\|_2^2 \leq \rho^{2k} \|x_0 - x_\ast\|_2^2. \label{eq: lin rate}
\end{align}
\end{theorem}
\begin{IEEEproof}
We proceed identically as in the proof of Theorems \ref{thm: weak thm} and \ref{thm: weak smooth thm}.  If \eqref{eq: strong LMI} is satisfied, then $e_k^\top Q_k e_k \leq 0$ for all $k$.  This is equivalent to $V_{k+1} -  \rho^2 V_k \leq 0$, by which \eqref{eq: lin rate} follows from Lemma \ref{lemma: lin lyapunov}.
\end{IEEEproof}

\section{Optimizing the Bound and Relaxation Parameter}

For each of the cases presented in the previous section we use the associated matrix inequality to optimize our bounds on the convergence rate.  In the case of non-strongly convex $f$ we provide analytic convergence rates, while for strongly convex $f$ we optimize our bounds numerically.

\subsubsection{Case 1: Non-strongly convex and non-smooth case}

\medskip
We now select algorithm parameters that satisfy the matrix inequality in \eqref{eq: weak LMI}.  In doing so, we arrive at a new and simple proof of the $O(1/k)$ convergence of DRS in the non-strongly convex and non-smooth case.

\begin{theorem} \label{theorem: weak nonsmooth}
If $f \in \mathcal F(0,\infty)$ and $g\in \mathcal F(0,\infty)$, then for any choice of $\lambda_k = \lambda \in (0,2)$ and $\alpha>0$, if we set $\sigma_k^{(1)} = \sigma_k^{(2)} = \sigma_k$, with
\begin{align}\label{opt sigma and theta nonsmooth}
\sigma_k:=2\lambda_k/\alpha, \quad \theta_k := \alpha^2\lambda_k(2-\lambda_k),
\end{align}
then $\sigma_k^{(1)}, \sigma_k^{(2)}, \alpha, \lambda_k,$ and $\theta_k$ satisfy the matrix inequality \eqref{eq: weak LMI}.
\end{theorem}
\begin{IEEEproof}
Making these substitutions gives $W_k^{(0)} + \sigma_k^{(1)}Q^{(1)} + \sigma_k^{(2)}Q^{(2)} = 0_n \preceq 0.$
\end{IEEEproof}


\begin{remark}
The convergence rate bound provided by the parameter choices in Theorem \ref{theorem: weak nonsmooth} guarantees convergence only if $\lim_{k \to \infty} \Theta_k = \infty$, which in this case means $\sum_{i=0}^\infty \lambda_k(2-\lambda_k) = \infty$.  This is consistent with the conditions on the relaxation parameter found in \cite{combettes2004solving}.
\end{remark}

\begin{remark} \label{rmk: lam 1}
After setting $\theta_k$ as in \eqref{opt sigma and theta nonsmooth}, we can maximize $\Theta_k$ by setting $\lambda_k = 1$ which results in $\Theta_k = \alpha^2 k$ and the following result
\begin{align} \label{eq: max result}
\min_{i=0,\ldots,k-1} \|\partial f(y_i) + \partial g(z_i)\|_2^2 \leq \frac{1}{\alpha^2 k}\|x_0-x_{\star}\|_2^2.
\end{align}
\end{remark}

\begin{remark}
Using the relation \eqref{eq: zy partial relation} we can rewrite \eqref{eq: max result} as
\begin{align}
\min_{i=0,\ldots,k-1} \|z_i - y_i\|_2^2 \leq \frac{1}{k}\|x_0-x_{\star}\|_2^2.
\end{align}
Thus we see that Theorem \ref{theorem: weak nonsmooth} also gives a $O(1/k)$ rate toward the iterates being a fixed point of the algorithm.
\end{remark}

\subsubsection{Case 2: Non-strongly convex and smooth case}

When $f \in \mathcal F(0,L_f), g\in \mathcal F(0,\infty)$ with $0<L_f<\infty$, we have the following result on feasibility of the matrix inequality \eqref{eq: weak smooth LMI}.

\begin{theorem} \label{thm: weak smooth thm param}
For any $\alpha > 0$ and $\lambda_k = \lambda $ with $0 < \lambda < 2$, if we set $\sigma_k^{(1)}=\sigma_k^{(2)} =\sigma_k$,
\begin{subequations}\label{opt sigma and theta smooth}
\begin{align}
\sigma_k &:=\frac{2\lambda_k}{\alpha}\left(\sqrt{\left(\frac{2-\lambda_k}{\alpha L_f}\right)^2 +1} - \left(\frac{2-\lambda_k}{\alpha L_f}\right)\right), \\
\theta_k &:= 2\lambda_k \alpha\! \left(1 \!+ \!\left(\frac{2-\lambda_k}{\alpha L_f}\right)\! - \! \sqrt{\left(\frac{2-\lambda_k}{\alpha L_f}\right)^2\! +\!1}\right). \label{eq: theta_1}
\end{align}
\end{subequations}
Then $\sigma_k^{(1)}, \sigma_k^{(2)}, \alpha, \lambda_k$, and $\theta_k$ satisfy the matrix inequality \eqref{eq: weak smooth LMI}.
\end{theorem}
\begin{IEEEproof}
This can be verified by substituting the expressions for $\sigma_k^{(1)}$, $\sigma_k^{(2)}$, and $\theta_k$ into the minors of $W_k^{(1)} + \sigma_k^{(1)}Q^{(2)} + \sigma_k^{(2)}Q^{(2)}$ and seeing that Sylverster's criterion is satisfied \cite{horn1990matrix}.
\end{IEEEproof}

\begin{remark}
Note that for a constant selection of $\theta_k = \theta$, $\Theta_k = \theta k$.  If \eqref{eq: weak smooth LMI} holds then
\begin{align}
	\min_{i=0,\ldots,k-1} [F(z_i) - F(z_\star)] \leq \dfrac{1}{\theta k} \|x_0 - x_\star\|_2^2.
	\end{align}
\end{remark}
\begin{remark}
For moderate values of $\alpha L_f$, we can take a second order Taylor expansion of the rightmost term in \eqref{eq: theta_1} and maximize the resulting expression with respect to $\lambda_k$.  This suggests that we should set $\lambda_k$ to
\begin{align}
\lambda_k = \frac{2}{3}\left(2 - \alpha L_f + \sqrt{1 - \alpha L_f + \alpha^2 L_f^2}\right).
\end{align}
\end{remark}

\subsubsection{Case 3: Strongly convex and smooth case}
\begin{figure}[ht]
\begin{center}
\vspace{-0.15in}
\centerline{\includegraphics[width=0.9\columnwidth]{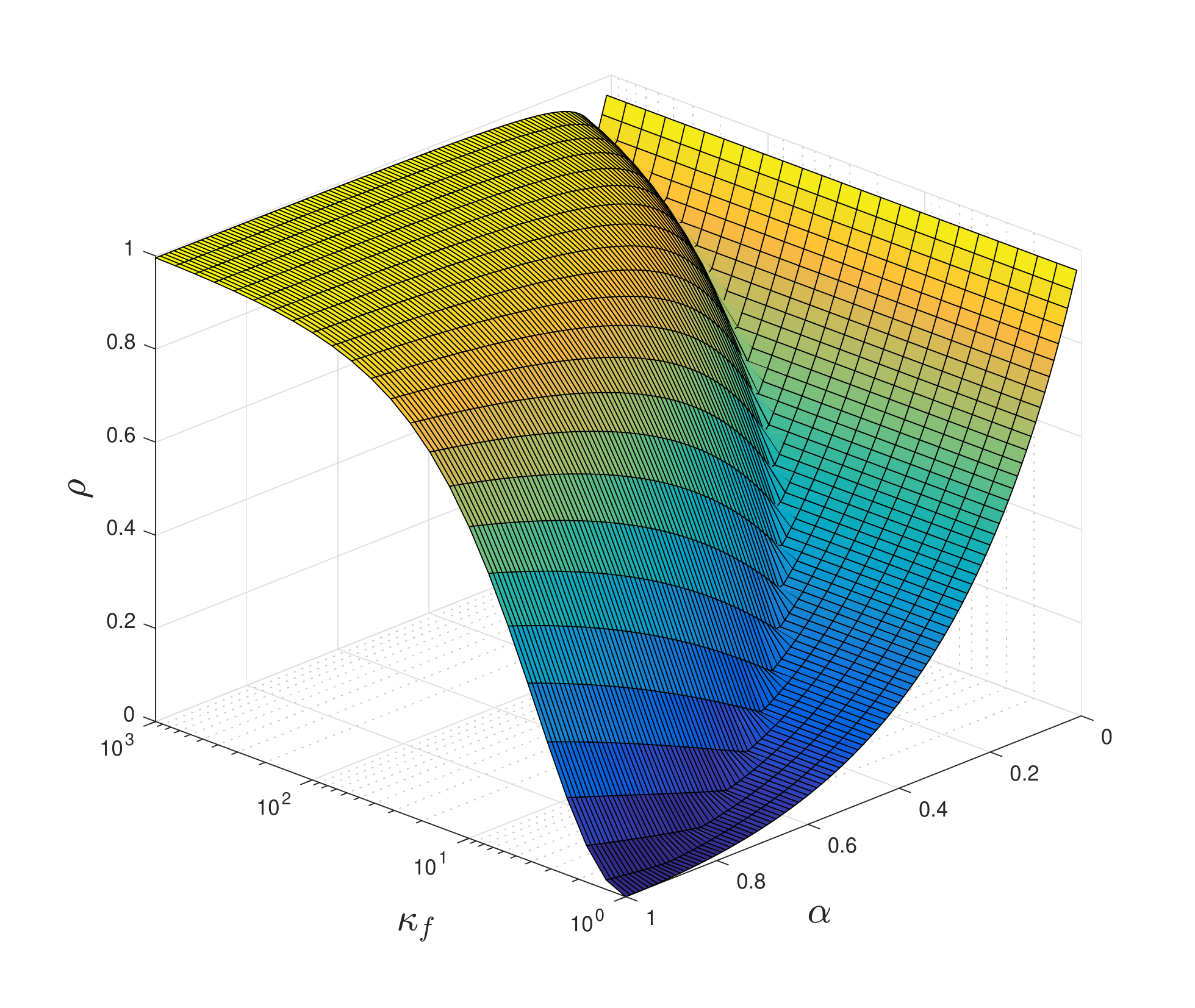}}
\caption{Optimal upper bound to linear convergence rate $\rho$ over $f \in \mathcal F(m_f,L_f)$, $g \in \mathcal F(0,\infty)$ as a function of step size $\alpha$ and condition number $\kappa_f = L_f/m_f$.}
\label{strong rho}
\end{center}
\end{figure}

When $f \in \mathcal F(m_f, L_f)$, $g\in\mathcal F(0,\infty)$, with $0 < m_f \leq L_f < \infty$, we can modify the matrix inequality \eqref{eq: strong LMI} to get a linear dependence on the relaxation parameter.  If we define $\Lambda_k := \begin{bmatrix} 0 & -\lambda_k & \lambda_k \end{bmatrix}^\top\otimes I_d$ and 
\begin{align}
M_k := (Q_k - \Lambda_k \Lambda_k^\top) +\sigma_k^{(1)} Q^{(1)} + \sigma_k^{(2)} Q^{(2)},
\end{align}
then \eqref{eq: strong LMI} is equivalent to
\begin{align} \label{eq: LMI V rewrite}
M_k - \Lambda_k [-1] \Lambda_k^\top \preceq 0.
\end{align}
As $\Lambda_k \Lambda_k^\top \succeq 0$, if \eqref{eq: LMI V rewrite} is satisfied then it must be the case that $M_k \preceq 0$.  We now recognize that $M_k- \Lambda_k [-1] \Lambda_k^\top$ is the Schur Complement of the bottom right entry in the matrix
\begin{align} \label{eq: Sig def}
\Sigma_k := \begin{bmatrix} M_k & \Lambda_k \\ \Lambda_k^\top & -1
\end{bmatrix}.
\end{align}
By the properties of the Schur complement \cite{boyd2004convex}, we can conclude that \eqref{eq: strong LMI} is satisfied if and only if $\Sigma_k \preceq 0$.  The advantage of using $\Sigma_k \preceq 0$ instead of \eqref{eq: strong LMI}, is that now both the convergence rate $ \rho^2$ and the relaxation parameter $\lambda_k$ appear linearly.  If we set $\sigma_k^{(1)} = \sigma^{(1)}$, $\sigma_k^{(2)} = \sigma^{(2)}$, and $\lambda_k = \lambda$ for all $k$, the corresponding bound in \eqref{eq: lin rate} can be optimized by solving the following SDP, where $\lambda$ is now a decision variable,
\begin{align}\label{eq: strong sdp}
\text{minimize} \quad\rho^2,\quad \text{subject to} \quad\Sigma_k \preceq 0,
\end{align} 
where the decision variables are $ \rho^2, \lambda > 0$ and $\sigma^{(1)},\sigma^{(2)} \geq 0$.  The optimal $\rho^2$ from solving this program over a range of step sizes $\alpha$ and condition numbers $\kappa_f = L_f/m_f$ is shown in Figure \ref{strong rho}.  We see that with increasing $\kappa_f$, the optimal choice of $\alpha$ decreases.  Across the range of values of $\kappa_f$ and $\alpha$, the SDP \eqref{eq: strong sdp} returns $\lambda=2$ as the optimal relaxation parameter.

\begin{remark}
We may repeat the same derivation for the inequalities \eqref{eq: weak LMI} and \eqref{eq: weak smooth LMI} to linearize their dependence on $\lambda_k$ as well.
\end{remark}

%

%
%

\section{Numerical Experiments}

We investigate how our theoretical results compare with the experimental performance of DRS in the three scenarios described above.  For the non-smooth and non-strongly convex case, we consider a basis pursuit problem (see \cite{boyd2011distributed}),
\begin{alignat}{2}
&\underset{x,z\in \R^n}{\text{minimize}} &&\quad \mathrm{1}_{\{y\in\R^n |Ay=b\}}(x) + \|z\|_1 \nonumber\\
&\text{subject to} &&\quad x - z = 0, \nonumber
\end{alignat}
with data $A \in \R^{300 \times 10000}$ and $b \in \mathbb{R}^{300}$.  We run DRS on the dual of this problem (ADMM) which is non-smooth and non-strongly convex.  We test the convergence over a range of values of $\lambda_k = \lambda$, including $\lambda^\star = 1$ (see Remark \ref{rmk: lam 1}), and fixed $\alpha = 1$.

\begin{figure*}[!htb]
\begin{center}
\centerline{\includegraphics[width=\linewidth]{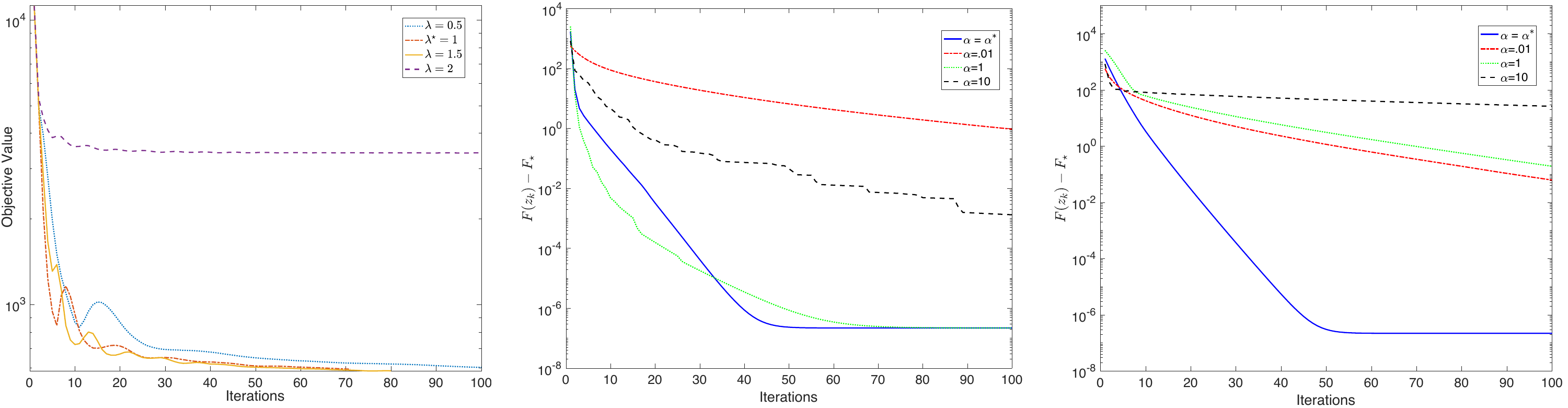}}
\centering\caption{From left to right: convergence of DRS on basis pursuit problem with $\alpha = 1$ and varying $\lambda$, convergence of DRS on LASSO problem with row-rank deficient $A$ with $\lambda = 1$ and varying $\alpha$, convergence of DRS on LASSO problem with full row-rank $A$ with fixed $\lambda = 2$ and varying $\alpha$.}
\label{fig: combo}
\end{center}
\end{figure*}

For the smooth cases we consider a LASSO problem,
\begin{align}
	\mathrm{minimize}_{x} \quad \frac{1}{2} \|Ax-b\|_2^2 + \gamma \|x\|_1.\nonumber
	\end{align}
with $A \in \R^{300 \times 200}$, $b \in \R^{300}$ and $\gamma = 0.1$.  For the non-strongly convex case we set $A$ to be rank deficient and plot the convergence of DRS over a range of $\alpha$ with $\lambda_k = 1$ fixed.  The value $\alpha^\star$ is found by performing a grid search over possible values and choosing that which maximizes the rate given by Theorem \ref{thm: weak smooth thm}.  For the strongly convex case $A$ is set to be full rank and plot the convergence of DRS over a range of $\alpha$ with $\lambda_k = 2$ fixed.  Again, $\alpha^\star$ is the value of $\alpha$ which gives the best rate as provided by Theorem \ref{thm: strong thm}.  The results are presented in Figure \ref{fig: combo}.

\section{Conclusion}

We presented a unified framework for deriving convergence bounds for DRS and parameter settings that optimize these bounds.  Our framework encompasses different assumptions on the smoothness and convexity of $f$.  We are able to give simple proofs of convergence and find optimal choices for the relaxation parameter by solving a small convex program for a fixed $\alpha$.  It is important to note that the parameter selections optimize our bounds in the sense of the best worst-case convergence rate over the entire class of objective functions with $f \in \mathcal F(m_f, L_f)$.  While there are scenarios where additional structure in the problem might make alternate parameter settings more effective, the settings we see here bound the worst-case performance agnostic of any additional problem structure.  For future work, this framework will be extended to encompass accelerated variants of DRS, as well as three or more operator splitting and multi-block ADMM.

%
%

\bibliographystyle{IEEEtran}
\bibliography{main}

%
\end{document}